\documentclass[12pt]{amsart}

\newcommand{\stopthm}{\hfill$\square$} 

\newcommand{\pa}{\partial}
\newcommand{\gh}{\widehat{g}} 
\newcommand{\Up}{\Upsilon}

\newcommand{\ep}{\epsilon}
\newcommand{\la}{\lambda}

\newcommand{\N}{{\mathbb N}}

\newcommand{\C}{{\mathbb C}}

\newcommand{\Ric}{\operatorname{Ric}}
\newcommand{\Res}{\operatorname{Res}}
\newcommand{\Pff}{\operatorname{Pff}}

\newcommand{\tr}{\operatorname{tr}}

\renewcommand{\Re}{\mathop{\rm Re}\nolimits}

\theoremstyle{plain}
\newtheorem{theorem}{Theorem}

\newtheorem{proposition}{Proposition}

\theoremstyle{definition}

\theoremstyle{remark}

\numberwithin{equation}{section}

\title{Holographic Formula for $Q$-Curvature}

\author{C. Robin Graham}
\address{Department of Mathematics, University of Washington,
Box 354350\\
Seattle, WA 98195 USA}
\email{robin@math.washington.edu}

\author{Andreas Juhl}
\address{Humboldt-Universit\"at, Institut f\"ur Mathematik, Unter den Linden, 
10099 Berlin}
\email{ajuhl@math.hu-berlin.de}

\begin{document}
\maketitle

\thispagestyle{empty}

\renewcommand{\thefootnote}{}
\footnotetext{The work of the first author was partially supported by NSF
grant DMS-0505701. The work of the second author was supported by 
SFB 647 ``Raum-Zeit-Materie'' of DFG.}

\section*{Introduction}\label{intro}

In this paper we give a formula for $Q$-curvature in even-dimensional 
conformal geometry. The $Q$-curvature was introduced by Tom Branson in
\cite{B} and has been the subject of much research. There are now a
number of characterizations of $Q$-curvature; see for example
\cite{GZ}, \cite{FG1}, \cite{GP}, \cite{FH}.  However, it has remained an
open problem to find an expression for $Q$-curvature which, for 
example, makes explicit the relation to the Pfaffian in the conformally  
flat case.  

\begin{theorem}\label{main}
The $Q$-curvature of a metric $g$ in even dimension $n$ is given by 
\begin{equation}\label{theformula}
2nc_{n/2}Q=nv^{(n)}+\sum_{k=1}^{n/2-1}(n-2k)p_{2k}^* v^{(n-2k)},   
\end{equation}
where $c_{n/2}=(-1)^{n/2}\left[ 2^n (n/2)!(n/2-1)!\right]^{-1}$.
\end{theorem}

Here the $v^{(2j)}$ are the coefficients appearing in the asymptotic
expansion of the volume form of a Poincar\'e metric for $g$, the
differential operators $p_{2k}$ are those which appear in the expansion of
a harmonic function for a Poincar\'e metric, and $p_{2k}^*$ denotes the
formal adjoint of $p_{2k}$.  These constructions are recalled 
in \S\ref{deriv} below.  We refer to the papers cited above and the
references therein for background about $Q$-curvature.  

Each of the operators $p_{2k}^*$ for $1\leq k\leq n/2-1$ can be factored as 
$p_{2k}^* = \delta q_k$, where $\delta$ denotes the divergence operator
with respect to $g$ and $q_k$ is a natural operator from functions to 
$1$-forms.  So the second term on the right hand side is the divergence of
a natural $1$-form.  In particular, integrating \eqref{theformula} over a
compact manifold recovers the result of \cite{GZ} that 
\begin{equation}\label{renormvol}
2c_{n/2}\int_M Q dv_g = \int_M v^{(n)} dv_g.
\end{equation}
This quantity is 
a global conformal invariant; the right hand side occurs as the coefficient  
of the log term in the 
renormalized volume expansion of a Poincar\'e metric (see \cite{G}). 

As we also discuss in \S\ref{deriv}, if $g$ is conformally flat then    
$$
v^{(n)} = (-2)^{-n/2} (n/2)!^{-1} \Pff,  
$$ 
where $\Pff$ denotes the Pfaffian of $g$. So in the conformally
flat case, Theorem~\ref{main} gives a decomposition of the
$Q$-curvature as a multiple of the Pfaffian and the divergence of a
natural 1-form.  A general result in invariant theory (\cite{BGP})
establishes the existence of such a decomposition, but does not produce a
specific realization.   

We refer to \eqref{theformula} as a holographic formula because its
ingredients come from the Poincar\'e metric, involving geometry in
$n+1$ dimensions.  Our proof is via the characterization of
$Q$-curvature presented in \cite{FG1} in terms of Poincar\'e metrics;
in some sense Theorem~\ref{main} is the result of making explicit the
characterization in \cite{FG1}. However, passing from the construction
in \cite{FG1} to \eqref{theformula} involves a non-obvious     
application of Green's identity. The transformation law of
$Q$-curvature under conformal change, probably its most fundamental
property, is not transparent from \eqref{theformula}, but it is from
the characterization in \cite{FG1}.  In \S\ref{relatedidentity}, we
derive another identity involving the $p_{2k}^*v^{(n-2k)}$ which is used in
\S \ref{families} and we discuss 
relations to the paper \cite{CQY}.  In 
\S\ref{families}, we describe the relation between holographic 
formulae for $Q$-curvature and the theory of conformally covariant   
families of differential operators of \cite{J}, and in particular explain
how this theory leads to the conjecture of a holographic formula for $Q$.         

We are grateful to the organizing committee of the 2007  Winter
School 'Geometry and Physics' at Srn\'i, particularly  to Vladimir
Sou\u{c}ek, for the invitation to this gathering, which made possible the
interaction leading to this paper.         

We dedicate this paper to the memory of Tom Branson.  His insights have led
to beautiful new mathematics and have greatly influenced our own respective
work.   

\section{Derivation}\label{deriv}

Let $g$ be a metric of signature $(p,q)$ on a manifold $M$ of even
dimension $n$. In this paper, by a Poincar\'e metric for $(M,g)$ we
will mean a metric $g_+$ on $M\times (0,a)$ for some $a>0$ of the form
\begin{equation}\label{formP}
g_+ = x^{-2}(dx^2 + g_x),
\end{equation}
where $g_x$ is a smooth 1-parameter family of metrics on $M$ satisfying 
$g_0=g$, such that $g_+$ is asymptotically Einstein in the sense that  
$\Ric(g_+) + ng_+ = O(x^{n-2})$ and 
$\tr_{g_+}(\Ric(g_+) +ng_+)=O(x^{n+2})$.  
Such a Poincar\'e metric always exists and $g_x$ is unique up addition
of a term of the form $x^n h_x$, where $h_x$ is a smooth 1-parameter family
of symmetric 2-tensors on $M$ satisfying $\tr_g(h_0) =0$ on $M$.  The
Taylor expansion of $g_x$ is even through order $n$ and the
derivatives $(\pa_x)^{2k}g_x|_{x=0}$ for $1\leq k\leq n/2-1$ and the trace   
$\tr_g((\pa_x)^ng_x|_{x=0})$ are determined inductively from the Einstein 
condition and are given by polynomial formulae in terms of $g$, its
inverse, and its curvature tensor and covariant derivatives thereof.  
See \cite{GH} for details. 

The first ingredient in our formula for $Q$-curvature consists of the
coefficients in the expansion of the volume form
\begin{equation}\label{volform}
dv_{g_+} = x^{-n-1}dv_{g_x}dx.
\end{equation}
Because the expansion of $g_x$ has only even terms through order $n$,
it follows that
\begin{equation}\label{volexp}
\begin{split}
dv_{g_x} &= \left(\frac{\det g_x}{\det g}\right)^{1/2}dv_{g}\\
&=(1+v^{(2)}x^2 + \cdots + v^{(n)}x^n +\cdots) dv_{g},
\end{split}
\end{equation}
where each of the $v^{(2k)}$ for $1\leq k \leq n/2$ is a smooth
function on $M$ expressible in terms of the curvature tensor of $g$
and its covariant derivatives. Set $v^{(0)}=1$.

The second ingredient in our formula is the family of differential
operators which appears in the expansion of a harmonic function for
the metric $g_+$. Given $f\in C^{\infty}(M)$, one can solve formally
the equation $\Delta_{g_+}u=O(x^n)$ for a smooth function $u$ such
that $u|_{x=0}=f$, and such a $u$ is uniquely determined modulo
$O(x^n)$.  The Taylor expansion of $u$ is even through order $n-2$ and
these Taylor coefficients are given by natural differential operators
in the metric $g$ applied to $f$ which are obtained inductively by
solving the equation $\Delta_{g_+}u=O(x^n)$ order by order.  See
\cite{GZ} for details.  We write the expansion of $u$ in the form
\begin{equation}\label{uexpansion}
u = f + p_2 f \, x^2 + \cdots + p_{n-2} f \,x^{n-2} + O(x^n); 
\end{equation}
then $p_{2k}$ has order $2k$ and its principal part is $(-1)^k
\dfrac{\Gamma(n/2-k)}{2^{2k}\,k!\,\Gamma(n/2)}\Delta^k$.  (Our
convention is $\Delta= -\nabla^i\nabla_i$.)  Set $p_0f=f$.

We remark that the volume coefficients $v^{(2k)}$ and the differential
operators $p_{2k}$ also arise in the context of an ambient metric
associated to $(M,[g])$. If an ambient metric is written in normal
form relative to $g$, then the same $v^{(2k)}$ are coefficients in the
expansion of its volume form, and the same operators $p_{2k}$ appear
in the expansion of a harmonic function homogeneous of degree 0 with
respect to the ambient metric.

Let $g_+$ be a Poincar\'e metric for $(M,g)$. In \cite{FG1} it is
shown that there is a unique solution $U \mod O(x^n)$ to
\begin{equation}\label{Uequation}
\Delta_{g_+} U = n +O(x^{n+1}\log x) 
\end{equation}
of the form
\begin{equation}\label{Uexp}
U=\log x + A + Bx^n \log x + O(x^n)\,,
\end{equation}
with
$$
A,B\in C^{\infty}(M\times [0,a))\,, \qquad A|_{x=0}=0\,. 
$$  
Also, $A \mod O(x^n)$ is even in $x$ and is formally determined by $g$, and   
\begin{equation}\label{BQ}
B|_{x=0}=-2c_{n/2}Q. 
\end{equation}
The proof of \eqref{BQ} presented in \cite{FG1} used results from
\cite{GZ} about the scattering matrix, so is restricted to positive
definite signature.  However, a purely formal proof was also indicated
in \cite{FG1}. Thus \eqref{BQ} holds in general signature.

\bigskip
\noindent
{\it Proof of Theorem~\ref{main}.} Let $g_+$ be a Poincar\'e metric
for $g$ and let $U$ be a solution of \eqref{Uequation} as described
above.
Let $f\in C^\infty(M)$ have compact support.  Let $u$ be a solution of 
$\Delta_{g_+}u=O(x^n)$ with $u|_{x=0}=f$; for  
definiteness we take $u$ to be given by \eqref{uexpansion} with the
$O(x^n)$ term set equal to $0$.  Let $0<\ep<x_0$ with $\ep$, $x_0$ small.    

Consider Green's identity   
\begin{equation}\label{green}
\int_{\ep<x<x_0}(U\Delta_{g_+}u - u\Delta_{g_+}U) \,dv_{g_+}
= \left(\int_{x=x_0}+\int_{x=\ep}\right)
\left(U\pa_{\nu}u - u\pa_{\nu}U\right)d\sigma,
\end{equation}
where $\nu$ denotes the inward normal and $d\sigma$ the induced volume
element on the boundary, relative to $g_+$.  Both sides have
asymptotic expansions as $\ep \rightarrow 0$; we calculate the
coefficient of $\log \ep$ in these expansions.

Using the form of the expansion of $U$ and the fact that
$\Delta_{g_+}u=O(x^n)$, one sees that the expansion of
$U\Delta_{g_+}u\,dv_{g_+}$ has no $x^{-1}$ term, so
$\int_{\ep<x<x_0}U\Delta_{g_+}u\,dv_{g_+}$ has no $\log\ep$ term.
Using \eqref{volform}, \eqref{volexp}, \eqref{uexpansion}, and
\eqref{Uequation}, one finds that the $\log\ep$ coefficient of
$-\int_{\ep<x<x_0}u\Delta_{g_+}U\,dv_{g_+}$ is
\begin{equation}\label{lhs}
n\sum_{k=0}^{n/2-1}\int_M v^{(n-2k)}p_{2k}f\,dv_g.
\end{equation}
On the right hand side of \eqref{green}, $\int_{x=x_0}$ is independent
of $\ep$, and
$$
\int_{x=\ep}\left(U\pa_{\nu}u - u\pa_{\nu}U\right)d\sigma
=\ep^{1-n}\int_{x=\ep}\left(U\pa_x u - u\pa_x U\right)dv_{g_\ep}.
$$ 
A $\log\ep$ term in the expansion of this quantity can arise only
from the $\log x$ or $x^n\log x$ terms in the expansion of $U$.
Substituting the expansions, one finds without difficulty that the
$\log\ep$ coefficient is
$$
\int_M \left(\sum_{k=1}^{n/2-1}2kv^{(n-2k)}p_{2k}f
-nBf\right)\,dv_g.
$$ 
Equating this to \eqref{lhs}, using \eqref{BQ}, and moving all
derivatives off $f$ gives the desired identity. \stopthm

\medskip
Since $\Delta_{g_+}1=0$, it follows that $p_{2k}1=0$ for $1\leq k\leq
n/2-1$.  Thus these $p_{2k}$ have no constant term, so
$p_{2k}^*=\delta q_{k}$ for some natural operator $q_{k}$ from
functions to 1-forms, where $\delta$ denotes the divergence with
respect to the metric $g$.  So in \eqref{theformula}, the second term
on the right hand side is the divergence of a natural $1$-form.  As
mentioned in the introduction, integration gives \eqref{renormvol}.
The proof of Theorem~\ref{main} presented above in the special case
$u=1$ is precisely the proof of \eqref{renormvol} presented in 
\cite{FG1}.

Theorem~\ref{main} provides an efficient way to calculate the $Q$
curvature. Solving for the beginning coefficients in the expansion of
the Poincar\'e metric and then expanding its volume form shows that
the first few of the $v^{(2k)}$ are given by:
\begin{align*}
\begin{aligned}
v^{(2)}&= -\frac12 J\\
v^{(4)}&= \frac18 (J^2-|P|^2)\\
v^{(6)}&= \frac{1}{48} \left(-\frac{2}{n-4}P^{ij}B_{ij} +3J|P|^2
-J^3 -2P^{ij}P_i{}^kP_{kj}\right) 
\end{aligned}
\end{align*}
where
\begin{align*}
\begin{aligned}
P_{ij}&=\frac{1}{n-2}\left(R_{ij}-\frac{R}{2(n-1)}g_{ij}\right)\\
J&=\frac{R}{2(n-1)}=P^i{}_i\\ 
B_{ij} &= P_{ij,k}{}^k -P_{ik,j}{}^k -P^{kl}W_{kijl}
\end{aligned}
\end{align*}
and $W_{ijkl}$ denotes the Weyl tensor. Similarly, one finds that the
operators $p_2$ and $p_4$ are given by:
\begin{align}
\begin{aligned}
-2(n-2)p_2&=\Delta\\
8(n-2)(n-4)p_4&= \Delta^2 +2J\Delta +2(n-2)P^{ij}\nabla_i\nabla_j
+(n-2)J,{}^i\nabla_i. 
\end{aligned}
\end{align}
For $n=2$, Theorem~\ref{main} states $Q=-2v^{(2)}=\frac12 R$.    
For $n=4$, substituting the above into Theorem~\ref{main} gives:
$$
Q=2(J^2 - |P|^2) + \Delta J,
$$
and for $n=6$:
\[
\begin{split}
Q & = 8P^{ij}B_{ij} +16 P^{ij}P_i{}^kP_{kj}-24J|P|^2+8J^3\\
& +\Delta^2J+4\Delta (J^2) +8(P^{ij}J_{,i})_{,j}-4\Delta (|P|^2). 
\end{split}
\]
In the formula for $n=6$, the first line is $(12c_{3})^{-1}6v^{(6)}$
and the second line is
$(12c_3)^{-1}\left(4p_2^*v^{(4)}+2p_4^*v^{(2)}\right)$.  Details of
these calculations will appear in \cite{J}.

The expansion of the Poincar\'e metric $g_+$ was identified explicitly
in the case that $g$ is conformally flat in \cite{SS}.  (Since we are only
interested in local considerations, by conformally flat we mean 
locally conformally flat.)  The two
dimensional case is somewhat anomalous in this regard, but the
identification of $Q$ curvature is trivial when $n=2$, so we assume
$n> 2$ for this discussion. The conclusion of \cite{SS} is that if
$g$ is conformally flat and $n>2$ (even or odd), then the expansion of
the Poincar\'e metric terminates at second order and
\begin{equation}\label{confflat}
(g_x)_{ij} = g_{ij} -P_{ij}x^2 +\frac14 P_{ik}P^k{}_j x^4.
\end{equation}
(The details of the computation are not given in \cite{SS}. Details 
will appear in \cite{FG2} and \cite{J}.) This easily yields  
\begin{proposition}\label{flatv}
If $g$ is conformally flat and $n> 2$, then     
$$
v^{(2k)}=\left\{ 
\begin{array}{ll} 
(-2)^{-k} \sigma_k(P) \qquad &  0\leq k\leq n\\
\qquad 0\qquad & n<k
\end{array} \right.
$$ 
where $\sigma_k(P)$ denotes the $k$-th elementary symmetric
function of the eigenvalues of the endomorphism $P_i{}^j$.
\end{proposition}
\begin{proof}
Write $g^{-1}P$ for $P_i{}^j$.  Then the $\sigma_k(P)$ are given by 
$$
\det(I+ g^{-1}P\,t) = \sum_{k=0}^n\sigma_k(P) t^k.  
$$
Equation \eqref{confflat} can be rewritten as 
$g^{-1}g_x = (I-\frac12 g^{-1}Px^2)^2$.  Taking the determinant and
comparing with \eqref{volexp} gives the result. 
\end{proof} 
\noindent
We remark that for $g$ conformally flat, $g_x$ given by \eqref{confflat} is
uniquely determined to all orders by the requirement that $g_+$ be
hyperbolic.     
So in this case the $v^{(2k)}$ are invariantly determined and given by
Proposition~\ref{flatv} for all $k\geq 0$ in all dimensions $n>2$.  

Returning to the even-dimensional case, we define the Pfaffian of the
metric $g$ by 
\begin{equation}\label{pff}  
2^n (n/2)!\,\Pff = (-1)^q \mu^{i_1\ldots i_n}\mu^{j_1\ldots j_n}  
R_{i_1i_2j_1j_2}\ldots R_{i_{n-1}i_nj_{n-1}j_n}, 
\end{equation}
where $\mu_{i_1\ldots i_n}=\sqrt{|\det(g)|}\,\ep_{i_1\ldots i_n}$  
is the volume form and $\ep_{i_1\ldots i_n}$ denotes the 
sign of the permutation. For a conformally flat metric, 
one has $R_{ijkl} =  2(P_{i[k}g_{l]j}-P_{j[k}g_{l]i})$. Using this in 
\eqref{pff} and simplifying gives  
$$
\Pff =  (n/2)!\,\sigma_{n/2}(P)  
$$ 
(see Proposition 8 of \cite{V} for details). Combining with
Proposition~\ref{flatv}, we obtain for conformally flat $g$:
$$
v^{(n)} = (-2)^{-n/2} (n/2)!^{-1} \Pff.   
$$
Hence in the conformally flat case, \eqref{theformula} specializes to   
$$
2Q = 2^{n/2} (n/2-1)!\Pff +   
(nc_{n/2})^{-1}\sum_{k=1}^{n/2-1}(n-2k)p_{2k}^*v^{(n-2k)},
$$
and again the second term on the right hand side is a formal divergence.

\section{A Related Identity}\label{relatedidentity}

In this section we derive another identity involving the 
$p_{2k}^*v^{(n-2k)}$.  It is in general  
impossible to choose the $O(x^n)$ term in \eqref{uexpansion} 
to make $\Delta_{g_+}u = O(x^n)$;
in fact $x^{-n}\Delta_{g_+}u|_{x=0}$ is independent of the $O(x^n)$ term
in \eqref{uexpansion} and is a conformally invariant operator of order $n$
applied to $f$, namely a multiple of the critical GJMS operator $P_n$.
Following 
\cite{GZ}, we consider the limiting behavior of the corresponding term in
the expansion of an eigenfunction for  
$\Delta_{g_+}$ as the eigenvalue tends to 0.  

Let $g_+$ be a Poincar\'e metric as above.
If $0\neq \lambda \in \C$ is near 0, then for $f\in C^{\infty}(M)$,
one can solve formally the equation 
$(\Delta_{g_+}-\la(n-\la))u_\la=O(x^{n+\la+1})$ for $u_\la$ of the form     
\begin{equation}\label{ula}
u_{\la} = x^\la\left( f +p_{2,\la}f\, x^2 + \cdots +p_{n,\la}f\, x^n
+O(x^{n+1})\right),
\end{equation}
where $p_{2k,\la}$ is a natural differential operator in the
metric $g$ of order $2k$ with principal part  
$(-1)^k\dfrac{\Gamma(n/2-k-\la)}{2^{2k}\,k!\,\Gamma(n/2-\la)}\Delta^k$ such
that $\dfrac{\Gamma(n/2-\la)}{\Gamma(n/2-k-\la)}p_{2k,\la}$ is polynomial
in $\la$.  Set $p_{0,\la}f=f$.   
The operators $p_{2k,\la}$ for $k<n/2$ extend analytically across $\la =0$  
and $p_{2k,0}=p_{2k}$ for such $k$, where $p_{2k}$ are the operators 
appearing in \eqref{uexpansion}.  But $p_{n,\la}$ has a simple pole at
$\la=0$ with residue a multiple of the critical GJMS operator 
$P_n$.  Now $P_n$ is self-adjoint, so it follows that  
$p_{n,\la}-p_{n,\la}^*$ is regular at $\la=0$.  We denote its value at
$\la=0$ by $p_n-p_n^*$, a natural operator of order at most  
$n-2$.  Our identity below involves the constant term $(p_n-p_n^*)1$.  
Note that since $P_n1=0$, both $p_{n,\la}1$ and 
$p_{n,\la}^*1$ are regular at $\la=0$.  We denote their values at $\la =0$
by $p_n1$ and $p_n^*1$; then $(p_n-p_n^*)1 = p_n1-p_n^*1$.  Moreover, 
(4.7), (4.13), (4.14) of \cite{GZ} show that 
\begin{equation}\label{pn1}
p_n1= -c_{n/2}Q.
\end{equation}
It is evident that $\int_M p_n1\, dv_g = \int_M p_n^*1\, dv_g$.  The next
proposition expresses the difference $p_n1-p_n^*1$ as a divergence.

\begin{proposition}\label{otheridentity}
\begin{equation}\label{theotheridentity}
n\left(p_n  -p_n^*\right)1 = \sum_{k=1}^{n/2-1}2k\,p_{2k}^*v^{(n-2k)}     
\end{equation}
\end{proposition}
\begin{proof}
Take $f\in C^\infty(M)$ to have compact support, let $0\neq \la$ be near 0,
and define $u_\la$ as in  
\eqref{ula} with the $O(x^{n+1})$ term taken to be 0.  Define $w_\la$ by
the corresponding expansion with $f=1$:
$$
w_{\la} = x^\la\left( 1 +p_{2,\la}1\, x^2 + \cdots +p_{n,\la}1\, x^n
\right). 
$$
As in the proof of Theorem~\ref{main}, consider Green's identity
\begin{equation}\label{green2}
\int_{\ep<x<x_0}(u_\la\Delta_{g_+}w_\la - w_\la\Delta_{g_+}u_\la) dv_{g_+} 
= \ep^{1-n}\int_{x=\ep}\left(u_\la\pa_x w_\la - w_\la\pa_x
u_\la\right)dv_{g_\ep} +c_{x_0}, 
\end{equation}
where $c_{x_0}$ is the constant (in $\ep$) arising from the boundary
integral over  
$x=x_0$.  Consider the coefficient of $\ep^{2\la}$ in the asymptotic
expansion of both sides.  The left hand side equals
$$
\int_{\ep<x<x_0}\left[u_\la\left(\Delta_{g_+}-\la(n-\la)\right)w_\la  
-w_\la\left(\Delta_{g_+}-\la(n-\la)\right)u_\la \right] 
 dv_{g_+}.
$$
Now $u_\la\left(\Delta_{g_+}-\la(n-\la)\right)w_\la\,dv_{g_+}$ and  
$w_\la\left(\Delta_{g_+}-\la(n-\la)\right)u_\la\,dv_{g_+}$ are
of the form $x^{2\la}\psi\, dxdv_g$ where $\psi$ is smooth up to $x=0$.    
It follows that the asymptotic expansion of the left hand side of
\eqref{green2} has no $\ep^{2\la}$ term.  Consequently the coefficient of
$\ep^{n+2\la}$ must vanish in the asymptotic expansion of 
$$
\int_{x=\ep}\left(u_\la x\pa_x w_\la - w_\la x\pa_x 
u_\la\right)dv_{g_\ep}.
$$
This is the same as the coefficient of $\ep^n$ in the expansion of   
\[
\begin{split}
\int_M&\left[\left(\sum_{k=0}^{n/2} p_{2k,\la}f\,\ep^{2k}\right)    
\left(\sum_{k=0}^{n/2} (2k+\la) p_{2k,\la}1\,\ep^{2k}\right)\right.\\  
&\qquad-\left.\left(\sum_{k=0}^{n/2} p_{2k,\la}1\,\ep^{2k}\right)  
\left(\sum_{k=0}^{n/2} (2k+\la) p_{2k,\la}f\,\ep^{2k}\right)
\right]\left(\sum_{k=0}^{n/2}v^{(2k)}\,\ep^{2k}\right)dv_g.
\end{split}
\]
Evaluation of the $\ep^n$ coefficient gives 
$$
\int_M \sum_{\substack{0\leq k,l,m\leq n/2\\k+l+m=n/2}}(2l-2k)(p_{2k,\la}f)
(p_{2l,\la}1) v^{(2m)}\, dv_g =0, 
$$
and then moving the derivatives off $f$ results in the pointwise identity  
\begin{equation}\label{lambdaidentity}
\sum_{\substack{0\leq k,l,m\leq n/2\\k+l+m=n/2}}(2l-2k)\,p_{2k,\la}^*
\left( (p_{2l,\la}1) v^{(2m)}\right) =0. 
\end{equation}
The limit as $\la\to 0$ exists of all $p_{2l,\la}1$ with $0\leq l\leq n/2$
and all $p_{2k,\la}^*$ with $0\leq k\leq n/2-1$.  Since $k=n/2$ forces
$l=m=0$, the operator $p_{n,\la}^*$ occurs only applied to 1.
Thus we may let $\la\to 0$ in \eqref{lambdaidentity}.  Using $p_{2l}1=0$ 
for $1\leq l\leq n/2-1$ results in
$$
np_n1 -\sum_{\substack{0\leq k,m\leq n/2\\k+m=n/2}} 2k\, p_{2k}^*v^{(2m)} 
=0. 
$$
Separating the $k=n/2$ term in the sum gives \eqref{theotheridentity}.  
\end{proof}
Proposition~\ref{otheridentity} may be combined with \eqref{theformula} and 
\eqref{pn1} to give other expressions for $Q$-curvature.  However,  
\eqref{theformula} seems the preferred form, as the other expressions all 
involve some nontrivial linear combination of $p_n1$ and $p_n^*1$.

We remark that the generalization of \eqref{lambdaidentity} obtained by  
replacing $p_{2l,\la}1$ by $p_{2l,\la}f$ remains true for arbitrary $f\in 
C^\infty(M)$.  This follows by the same argument, taking $w_\la$ to be 
given by the asymptotic expansion of the same form but with arbitrary 
leading coefficient.  

We conclude this section with some observations concerning relations to the
paper \cite{CQY}:  
\begin{enumerate}
\item  Recall that Theorem~\ref{main} was proven by consideration of the 
$\log\ep$ term in \eqref{green}, generalizing the proof of
\eqref{renormvol} in \cite{FG1} where $u=1$.  In \cite{CQY}, it was shown
that for a  
global conformally compact Einstein metric $g_+$, consideration
of the constant term in 
$$
\int_{x>\ep}\Delta_{g_+}U \,dv_{g_+}
=\int_{x=\ep} \pa_{\nu}U \,d\sigma
$$
for $U$ a global  
solution of $\Delta_{g_+}U =n$ gives a formula for the 
renormalized volume $V(g_+,g)$ of $g_+$ relative to a metric $g$ in the
conformal infinity of $g_+$.  In our notation this formula reads 
\begin{equation}\label{cqyformula}
V(g_+,g)=-\int_M \frac{d}{ds}\Big|_{s=n}(S(s)1)\,dv_g    
+\frac{1}{n}\int_M \sum_{k=1}^{n/2}2k\,\dot{p}_{2k}^* v^{(n-2k)}\,dv_g,
\end{equation}
where $\dot{p}_{2k}=\frac{d}{d\la}|_{\la =0}p_{2k,\la}$ (which exists for
$k=n/2$ when applied to 1) and $S(s)$ denotes  
the scattering operator relative to $g$.  The operators $\dot{p}_{2k}$
arise in this context because the coefficient of $x^{2k}$ in the expansion
of $U$ is $\dot{p}_{2k}1$ for $1\leq k\leq n/2-1$, and the coefficient of
$x^n$ involves $\dot{p}_{n}1$.  Likewise, consideration  
of the constant term in 
$$
\int_{x>\ep} u\Delta_{g_+}U \,dv_{g_+}
=\int_{x=\ep} \left( u\pa_{\nu}U-U\pa_{\nu}u \right)d\sigma 
$$
for harmonic $u$ gives an analogous formula for the finite part of  
$\int_{x>\ep}u \,dv_{g_+}$ in terms of boundary data.  
\item There is an 
analogue of Proposition~\ref{otheridentity} involving the 
$\dot{p}_{2k}^* v^{(n-2k)}$.  Differentiating 
\eqref{lambdaidentity} with respect   
to $\la$ at $\la=0$ and rearranging gives the identity  
$$
\sum_{k=1}^{n/2}2k \left(\dot{p}_{2k}^* 
v^{(n-2k)}-(\dot{p}_{2k}1)v^{(n-2k)}\right)
=\sum_{k=2}^{n/2}\sum_{l=1}^{k-1}(4l-2k)p_{2k-2l}^* 
\left((\dot{p}_{2l}1)v^{(n-2k)}\right)
$$  
which expresses the left hand side as a divergence.  
\item In \cite{CQY} it was also shown that under an infinitesimal conformal 
change, the scattering term 
$$
\mathcal{S}(g_+,g)\equiv\int_M \frac{d}{ds}\Big|_{s=n}(S(s)1)dv_g   
$$
satisfies
$$
\frac{d}{d\alpha}\Big|_{\alpha =0}\, \mathcal{S}(g_+,e^{2\alpha\Up}g)   
=-2c_{n/2}\int_M \Up Q \,dv_g.
$$
Comparing with 
$$
\frac{d}{d\alpha}\Big|_{\alpha =0}\, V(g_+,e^{2\alpha\Up}g)  = 
\int_M \Up v^{(n)} \,dv_g
$$
(see \cite{G}) and using \eqref{cqyformula} and Theorem~\ref{main},  
one deduces the curious conclusion that the infinitesimal conformal
variation of  
$$
\int_M \sum_{k=1}^{n/2}2k\,\dot{p}_{2k}^* v^{(n-2k)}\,dv_g
$$
is
$$
-\int_M \Up \sum_{k=1}^{n/2-1}(n-2k)p_{2k}^* v^{(n-2k)} \, dv_g.
$$
This statement involves the conformal variation only of local expressions.      
For $n=2$ this is the statement of conformal invariance of $\int_M R\, 
dv_g$, 
while for $n=4$ it is the assertion that the infinitesimal conformal
variation of  $\int_M J^2\, dv_g$ is $2\int_M \Up \Delta J\, dv_g$.
\end{enumerate}

\section{$Q$-curvature and families of conformally covariant differential operators}
\label{families}

In \cite{J} one of the authors initiated a theory of one-parameter
families of natural conformally covariant local operators
\begin{equation}\label{family}
D_N(X,M;h;\lambda): C^\infty(X) \to C^\infty(M), \quad N \ge 0
\end{equation}
of order $N$ associated to a Riemannian manifold $(X,h)$ and a hypersurface 
$i:M \to X$, depending rationally on the parameter $\la\in \C$.  For 
such  
a family the conformal weights which describe the covariance of the
family are coupled to the family parameter in the sense that
\begin{equation}\label{CI}
e^{-(\lambda-N)\omega} D_N(X,M;\widehat{h};\lambda) e^{\lambda \omega} =
D_N(X,M;h;\lambda), \quad \widehat{h} = e^{2\omega} h
\end{equation}
for all $\omega \in C^\infty(X)$ (near $M$). 

Two families are defined in \cite{J}: 
one via a residue construction which has its origin in an extension
problem for automorphic functions of Kleinian groups through their
limit set (\cite{J2}, chapter 8), and the other via a tractor
construction.
Whereas the tractor family depends on the choice of a metric $h$   
on $X$, the residue family depends on the choice of an  
asymptotically hyperbolic metric $h_+$ and a defining function $x$, to
which is associated the metric $h=x^2h_+$.    

Fix an asymptotically hyperbolic metric $h_+$ on one side $X_+$ of $X$
in $M$ and choose a defining function $x$ for $M$ with $x>0$ in $X_+$.
Set $h=x^2h_+$. To an eigenfunction $u$ on $X_+$ satisfying
$$
\Delta_{h_+}u = \mu(n-\mu)u, \qquad  \Re\mu=n/2, \quad \mu \neq n/2 
$$
is associated the family
$$
\langle T_u(\zeta,x),\varphi\rangle \equiv 
\int_{X_+} x^\zeta \, u \, \varphi \,dv_{h},\quad \varphi\in C_c^\infty(X) 
$$ 
of distributions on $X$. The integral converges for
$\Re\zeta>-n/2-1$ and the existence of a formal asymptotic expansion
$$
u\sim \sum_{j\geq 0}x^{\mu +j}a_j(\mu) + 
\sum_{j\geq 0}x^{n-\mu +j}b_j(\mu),\quad x\to 0 
$$
with $a_j(\mu), b_j(\mu)\in C^\infty(M)$ 
implies the existence of a meromorphic continuation of $T_u(\zeta,x)$ to
$\C$ with simple poles in the ladders 
$$
-\mu -1-\N_0,\qquad -(n-\mu)-1-\N_0.
$$
For $N\in \N_0$, its residue at $\zeta = -\mu - 1 -N$ has the form   
$$
\varphi\mapsto \int_M a_0 \delta_N(h;\mu+N-n)(\varphi)dv_{i^*h}, 
$$
where 
$$
\delta_N(h;\lambda): C^\infty(X)\to C^\infty(M)
$$ 
is a family of differential operators of order $N$ depending 
rationally on $\lambda \in \C$. If $\widehat{x}  
= e^{\omega}x$ with $\omega \in C^\infty(X)$, then
$\widehat{h}=e^{2\omega}h$ and it is easily checked that
$\delta_N(h;\lambda)$ satisfies \eqref{CI}.  (The family
$\delta_N(h;\lambda)$ should more correctly be regarded as determined
by $x$ and $h_+$, but we use this notation nonetheless.)

If $g$ is a metric on $M$, then we can take $h_+=g_+$ to be a
Poincar\'e metric for $g$ on $X_+=M\times (0,a)$ and $x$ to be the
coordinate in the second factor, so that $h=dx^2 +g_x$. Then
(assuming $N \leq n$ if $n$ is even), the family $\delta_N(h;\lambda)$
depends only on the initial metric $g$. The residue can be evaluated
explicitly and for even orders $N=2L$ one obtains 
\begin{equation}\label{res-critical}
\delta_{2L}(h;\mu+2L-n)=  
\sum_{k=0}^{L} \frac{1}{(2L-2k)!} \left(
\sum_{l=0}^k p_{2l,\mu}^*\circ v^{(2k-2l)} \right)\circ 
i^* \partial_x^{2L-2k},  
\end{equation}
where the $p_{2l,\mu}$ are the operators appearing in \eqref{ula}  
and the coefficients $v^{(2j)}$ are used as multiplication operators. 
The corresponding residue family is defined by 
\begin{equation}\label{resdef}
D_{2L}^{res}(g;\lambda)=2^{2L} L!  \frac{\Gamma(-n/2-\lambda+2L)}
{\Gamma(-n/2-\lambda +L)}\, \delta_{2L}(h;\lambda);  
\end{equation}
the normalizing factor makes $D_{2L}^{res}(g;\lambda)$ polynomial in 
$\lambda$. 
We are interested in the critical case $2L=n$ for $n$ even.  
Using
$$
\Res_0(p_{n,\lambda}) = - c_{n/2} P_n
$$ 
from \cite{GZ}, we see that   
\begin{equation}\label{res-value}
D_n^{res}(g;0) = (-1)^{n/2}P_n(g) i^*.
\end{equation}
Direct evaluation from \eqref{res-critical}, \eqref{resdef} gives    
$$
\dot{D}^{res}_n(g;0)1 = -(-1)^{n/2} c_{n/2}^{-1} \left(
p_n^* 1 + \sum_{k=0}^{n/2-1} p_{2k}^*v^{(n-2k)} \right), 
$$
where the dot refers to the derivative in $\la$.  

Suppose now that $g$ is transformed conformally:
$\widehat{g}=e^{2\Up}g$ with $\Up \in C^\infty(M)$.  By the construction of
the normal form in \S 5 of \cite{GL}, the 
Poincar\'e metrics $g_+$ and $\widehat{g}_+$ are related by $\Phi^*
\widehat{g}_+ = g_+$ for a diffeomorphism $\Phi$ which restricts to
the identity on $M$ and for which the function $\Phi^*(x)/x$ restricts to   
$e^{\Up}$. Using this the residue construction easily implies 
\begin{equation}\label{CIP}
e^{-(\lambda-n)\Up} D_n^{res}(\gh;\lambda) =
D_n^{res}(g;\lambda) \left(\Phi^*(x)/x\right)^{-\lambda} \Phi^*.  
\end{equation}
Applying \eqref{CIP} to the function 1, differentiating at $\la =0$, and
using \eqref{res-value} and $P_n1=0$ gives
$$
e^{n\Up}\dot{D}^{res}_n(\gh;0)1= \dot{D}^{res}_n(g;0)1 -
(-1)^{n/2}P_{n/2}\Up.  
$$
This proves that the curvature quantity 
$$
-(-1)^{n/2}\dot{D}^{res}_n(g;0)1=c_{n/2}^{-1} \left(
p_n^* 1 + \sum_{k=0}^{n/2-1} p_{2k}^*v^{(n-2k)}\right)
$$ 
satisfies the same transformation law as the $Q$-curvature. It is
natural to conjecture that it equals the $Q$-curvature.  Indeed, this
follows from \eqref{theformula}, \eqref{pn1}, and
\eqref{theotheridentity}:
\[
\begin{split}
p_n^* 1 + \sum_{k=0}^{n/2-1} p_{2k}^*v^{(n-2k)}=
 p_n1 &+
\left( p_n^* 1-p_n1 + \frac{1}{n}\sum_{k=1}^{n/2-1} 2k
  p_{2k}^*v^{(n-2k)}\right)\\ 
&
+\left( \sum_{k=0}^{n/2-1}
  p_{2k}^*v^{(n-2k)}-\frac{1}{n}\sum_{k=1}^{n/2-1}2k
  p_{2k}^*v^{(n-2k)}\right). 
\end{split}
\]
\noindent
The first term is $-c_{n/2}Q$ by \eqref{pn1},
the second term is $0$ by \eqref{theotheridentity}, and 
the last term is $2c_{n/2}Q$ by \eqref{theformula}.  

The relation $\dot{D}^{res}_n(g;0)1=(-1)^{n/2+1}Q$ and 
\eqref{res-value} show that both the critical GJMS operator 
$P_n$ and the $Q$-curvature are contained in the one 
object $D_n^{res}(g;\lambda)$.  In that respect, $D_n^{res}(g;\lambda)$ 
resembles the scattering operator in \cite{GZ}.  However, the family 
$D_n^{res}(g;\lambda)$ is local and all operators in the family have  
order $n$.

\end{document}